\documentclass[12pt,reqno,a4paper]{amsart}
\usepackage{mathrsfs}
\usepackage{extsizes}
\usepackage{blindtext}
\usepackage{amsmath,amsfonts,amsthm}
\usepackage{amssymb,latexsym}
\usepackage[colorlinks]{hyperref}
\usepackage{graphicx}
\usepackage[all]{xy}

\usepackage{layout}
\usepackage{bbm}
\usepackage[pagewise]{lineno}
\usepackage{tikz}

\usetikzlibrary{arrows.meta,
	decorations.markings,
	positioning}

\theoremstyle{plain}
\newtheorem{theorem}{Theorem}[section]

\theoremstyle{definition}

\newtheorem{definition}[theorem]{Definition}

\theoremstyle{remark}

\numberwithin{equation}{section}

\newcommand{\LL}{\mathscr{L}}
\newcommand{\K}{\mathscr{K}}
\newcommand{\la}{\langle}
\newcommand{\ra}{\rangle}

\begin{document}
 	\title[On Simplicity of Cuntz Algebra and its generalizations]{On Simplicity of Cuntz Algebra and its generalizations}
 	
\author{Massoud Amini}
\author{Mahdi Moosazadeh}
\address{Department of Mathematics,
	Faculty of Mathematical Sciences,
	Tarbiat Modares University, Tehran 14115-134, Iran}
\curraddr{}
\email{mamini@modares.ac.ir, mahdimoosazadeh7@gmail.com}
\thanks{}

\thanks{}

 \begin{abstract}
Cuntz algebra $\mathcal O_2$ is the universal $C^*$-algebra generated by two isometries $s_1, s_2$ satisfying $s_1s_1^*+s_2s_2^*=1$. This is separable, simple, infinite $C^*$-algebra containing a copy of any nuclear $C^*$-algebra. The $C^*$-algebra $\mathcal O_2$ plays a central role in the modern theory of $C^*$-algebras and appears in many substantial statements, including a formulation of the celebrated Uniform Coefficient Theorem (UCT). There are several extensions of this notion, including Cuntz algebra $\mathcal O_n$, Cuntz-Krieger algebra $\mathcal O_A$ for a matrix $A$, Cuntz-Pimsner algebra $\mathcal O_X$ and its relaxation by Katsura for a $C^*$-correspondence $X$, and Cuntz-Nica-Pimsner algebra $\mathcal {NO}_X$, for a product system $X$. We give an overview of the construction of these classes of $C^*$-algebras with a focus on conditions ensuring their simplicity, which is needed in the Elliott Classification Program, as it stands now. The results we present are now part of the literature, except our discussion on the sufficient conditions for simplicity of the reduced Cuntz-Nica-Pimsner algebra $\mathcal{NO}^r_X$, which is known to experts, but might happen to be new for some of our audiences. 
\end{abstract}	


\date{}

\dedicatory{}

\maketitle

\section{Introduction} 	

The Cuntz algebra ${\mathcal {O}}_{n}$, introduced in the late 70's by Joachim Cuntz \cite{c}, is the universal $C^*$-algebra generated by 
$n$ isometries $s_1,\cdots s_n$ of an infinite-dimensional Hilbert space 
satisfying $\sum_{i=1}^{n} s_is_i^*=1$. These are the first concrete examples of a separable infinite simple $C^*$-algebra and any such algebra contains a subalgebra that has some $\mathcal {O}_{n}$ as a quotient. In the early 80's, Joachim Cuntz and Wolfgang  Krieger extended the class of Cuntz algebras to a class of simple $C^*$-algebras generated by partial isometries,  known as the Cuntz-Krieger algebras \cite{ck}. These are purely infinite, but not necessarily simple. For an $n\times n$ matrix $A$ with entries $0$ or $1$ and non zero rows and columns, the Cuntz-Krieger algebra $\mathcal O_A$ is the universal $C^*$-algebra generated by 
$n$ partial isometries $s_1,\cdots s_n$ of an infinite-dimensional Hilbert space 
satisfying $\sum_{i=1}^{n} s_is_i^*=1$ and $s_i^*s_i=\sum_{j=1}^{n} A_{i,j}s_js_j^*$, for $i=1,\cdots, n$. The Cuntz-Krieger algebra $\mathcal O_A$ is simple (and thereby  independent of the choice of
generators) if and only if $A$ is irreducible (i.e., for any pre choice of a row and column indicies some power of $A$ has a nonzero entry at that row and column) and not a permutation matrix (i.e., no power of $A$ is the identity).  
In the late 90's, Mihai Pimsner \cite{pim} developed a far reaching generalization of previous constructions
by associating to each $C^*$-correspondence $X$ over a $C^*$-algebra $A$ two $C^*$-algebras $\mathcal T_X$ and $\mathcal O_X$, respectively extending Toeplitz and Cuntz-Krieger algebras. The Cuntz-Pimsner algebras $\mathcal O_X$ also generalize crossed products by $\mathbb Z$. For a $\mathbb Z$-$C^*$-algebra, if we endow $X=A$ with the canonical $A$-valued inner product and right $A$-module structure $a\cdot b:=a(1\cdot b)$, then the Cuntz-Pimsner algebra $\mathcal O_X$ is isomorphic to the
full crossed product $A\rtimes\mathbb  Z$. J\"urgen Schweizer \cite[Theorem 3.9]{s} observed that if $A$ is unital and $X$ is full $A$-correspondence with injective left
action, then $\mathbb O_X$ is simple if and only $X$ is aperiodic (i.e., the $n$-fold tensor $X^{\otimes\ n}$ is unitarily equivalent to the identity correspondence over $A$ only if $n=0$) and $A$ has no nontrivial
invariant ideal. Later, Menev\c se Ery\"uzl\"u and Mark Tomforde \cite[Theorem 4.3]{et}  showed that for $\mathcal O_X$ being simple, one could replace the Schweizer’s aperiodicity condition with the so called $(S)$ condition, that is, for each positive element $a\in A_+$, $\varepsilon>0$, and $n\geq 1$, there is $m > n$ and a non returning unit vector $\zeta\in X^{\otimes m}$  (i.e., as elements in $\mathcal O_X$ via universal representation, $\zeta^*\xi\zeta=0$, for  $\xi\in X^{\otimes n}$) with $\|\langle a\zeta,\zeta \rangle\|>\|a\|-\varepsilon$. 

In general, a right Hilbert $A$-bimodules behaves like a  generalized endomorphisms of $A$ and thereby $\mathcal T_X$ and $\mathcal O_X$ looks respectively like $A\rtimes\mathbb N$ and $A\rtimes\mathbb Z$ \cite{clsv}. Inspired by the graph $C^*$-algebras, Takeshi Katsura
\cite{k} extended Pimsner approach to modules $X$ with not necessarily isometric left actions, and defined the most general version of the $C^*$-algebra $\mathcal O_X$ for an $A$-correspondence $X$, and observed that the Cuntz-Pimsner representation of $X$ generates an isomorphic copy of $\mathcal O_X$ when it is injective
and respects the gauge action.
In another direction, adapting a notion product systems over
of Hilbert spaces on the interval $(0,\infty)$ due to William Arveson \cite{a},  Neal J. Fowler  introduced a notion of $C^*$-algebras associated to product systems of Hilbert $C^*$-bimodules \cite{f2} (see also, \cite{d,f2} for the same notion on a general  semigroup). A product system over a semigroup
$P$ of Hilbert $A$–bimodules is nothing but an action of $P$ on $A$ by generalized
endomorphisms. Alexandru Nica \cite{n} had already introduced certain Toeplitz type algebras for semigroups $P$ sitting inside a group $G$ inducing a quasi-lattice order.  and Fowler was able to associate a Toeplitz  $C^*$-algebra $\mathcal T_X^{\rm cov} $ and a Cuntz-Pimsner $C^*$-algebra $\mathcal O_X$ to a so called compactly aligned
product system $X$ over a quasi-lattice ordered group $(G, P)$, with $\mathcal O_X$ no longer a quotient of $\mathcal T^{\rm cov}_X$ and the canonical morphism from $A$ to $\mathcal O_X$ no longer injective in general,  and show  that it resembles twisted crossed product $A\rtimes_\sigma P$ for twist $\sigma$ coming from $X$.
Aidan Sims and Trent Yeend \cite{sy} introduced the
Cuntz-Nica-Pimsner algebra $\mathcal{NO}_X$ of a product system $X$ over a quasi-lattice ordered group $(G,P)$, which happens to be a quotient of $\mathcal T^{\rm cov}_X$ and the canonical representation of $X$ on $\mathcal{NO}_X$ is isometric under 
certain mild conditions. Finally, Toke M. Carlsen, Nadia S. Larsen, Aidan Sims and Sean T. Vittadello \cite{clsv} managed to describe  the core of $\mathcal{NO}_X$  as well as  the canonical coaction of $G$ on $\mathcal{NO}_X$ and prove a gauge-invariance uniqueness theorem for  $\mathcal{NO}_X$, where Katsura’s gauge action of $\mathbb T$ (that is, a coaction of $\mathbb Z$) is replaced by a coaction of $G$.
Moreover, for a large class of product systems $X$, they construct a reduced version $\mathcal {NO}^r_X$,  generated by an injective Nica covariant Toeplitz representation of $X$, which admits a coaction of
$G$ compatible with that of the gauge group on $\mathcal T^{\rm cov}_X$, satisfying  the expected couniversal property (like that of Exel \cite{e} for Fell bundles and Katsura \cite{k} for $C^*$-algebras of correspondences). This is the same as the full Cuntz-Nica-Pimsner algebra $\mathcal{NO}_X$ under certain amenability type condition.

As far as we know, there are few (if any) surveys on simplicity of these important classes of $C^*$-algebras. In particular, we could not trace any article directly discussing simplicity of  Cuntz-Nica-Pimsner algebras, though as explained in the last section, necessary and sufficient conditions are essentially known.

\section{From Cuntz Algebra to Higher rank graph C$^*$-algebra}

The notion of a {\it universal} C$^*$-algebra of a set of generators and relations is needed to define most of the later constructions. Let us start by recalling this construction due to Bruce Blackadar \cite{b}. First thing to note is that such a universal $C^*$-algebra does not exist for any choice of relations, as arbitrary set of relations may not realizable by operators on Hilbert spaces. This suggests that one should start with a notion of representation.

	Let $\mathscr G$  be a set of generators and  $\mathscr R$ be a set of relations between generators (which are meaningful for operators on Hilbert spaces). A representation $\pi : \mathscr G \rightarrow \mathbb B(\mathscr{H})$ of $(\mathscr G , \mathscr R)$ on a Hilbert space $\mathscr{H}$ is a map such that $\pi(\mathscr G)$ satisfies the same relations as $\mathscr R$ in $\mathbb B(\mathscr H)$.

Let $\mathscr A$ be the free $*$-algebra on  $\mathscr G$, then a representation $\pi$ of $(\mathscr G , \mathscr R)$ on $\mathscr{H}$ algebraically extends to a $*$-representation of $\mathscr A$ on $\mathscr{H}$. For $x \in \mathscr A$, let us define
\[
	\| x \| := \sup \{ \| \pi(x) \|:\ \pi \text{ is a representation of } (\mathscr G , \mathscr R) \}.
\]

Now if this supremum is finite, for all $x \in \mathscr{G}$, it  defines a C$^*$-seminorm on $\mathscr A$. The completion  $C^*(\mathscr G, \mathscr R)$ of $\mathscr{A}$ is then a $C^*$-algebra, called the universal  C$^*$-algebra of $\mathscr{G}$ subjected to relations $\mathscr{R}$. A typical instance of relations leading to a finite seminorm is one of the form
$$\|p(x_{i_1},\cdots, x_{i_n}, x^*_{i_1},\cdots, x^*_{i_n})\|\leq \varepsilon,$$ 
for a polynomial $p$ of $2n$ noncommuting variables with complex coefficients. 

For a finite set of generators, a finite set of such relations, plus a finite set of finite algebraic relations, the universal $C^*$-algebra is known to exist, though the case of infinitely many generators and relations is not ruled out.
    If the universal C$^*$-algebra $C^*(\mathscr G, \mathscr R)$ exists, it has a universal property: for any C$^*$-algebra $\mathcal A$ with a subset $X$ in a one-to-one correspondence with $\mathscr G$, satisfying relations as those of $\mathscr R$, there is a surjective $*$-homomorphism from $C^*(\mathscr G, \mathscr R)$ onto the $C^*$-subalgebra of $\mathcal A$ generated by $X$.

Now we are ready to describe the classes of $C^*$-algebras considered in this survey. Generally, these are constructed as universal C$^*$-algebras for some given sets of generators and relations. We discuss the simplicity results in each case separately.

\subsection{Cuntz Algebras}
We start with  Cuntz algebras, introduced in \cite{c} by Joachim Cuntz, which are unital, simple, nuclear, and purely infinite. 

\begin{definition}
	Let $\mathscr H$ be a separable Hilbert space, $n \geq 2$ and $s_1, \dots s_n$ be isometries in $B(\mathscr H)$, that is, $s_i^*s_i = id_{\mathscr H} =: 1$, for $i=1,\cdots, n$.  The universal $C^*$-algebra generated by $\{s_i\}_{i = 1}^{n}$ subjected to a single relation $\sum_{i = 1}^{n} s_is_i^* = 1$, is called the Cuntz algebra $\mathcal O_n$. It follows from definition that $s_i^*s_j = 0$, if $i \neq j$, that is, the range projections of the distinct generators are orthogonal.
	
\end{definition}

For the Cuntz algebra $\mathcal O_n$ with generators $s_1, \cdots s_n$, let $W$ be the set of all $k$-tuples $(i_1,\cdots, i_k)$, for $k\geq 1$, with $i_j \in \{ 1, \cdots, n\}$, for $j=1,\cdots k$, let $\ell : W \rightarrow \mathbb{N}$ to the corresponding length function. To each $(i_1, \dots i_k) =: \mu \in W$, associate an isometry $s_\mu := s_{1_1} s_{i_2} \dots s_{i_k}$. The set $E_k := \{ s_\mu s_\nu^*:\ \ell(\mu) = \ell(\nu) \}$ forms a matrix unit system for $\mathbb M_{n^k}(\mathbb{C})$, which then spans $F_k := \mathbb M_{n^k}(\mathbb{C})$. Since $\{E_k\}$ is an increasing sequence, $F = \bigcup_{k \in K} F_{k}$ is a UHF-algebra $\mathbb M_{n^\infty}$ of type $n^\infty$. Now for the range projection $p_1:=s_1s^*_1$, the cutdown $\rho_1$ by $p_1$ gives rise to an action of $\mathbb N$ on $\mathbb M_{n^\infty}$ by endomorphisms, and the Cuntz algebra $\mathcal O_n$ is isomorphic to the corresponding semigroup crossed product $\mathbb M_{n^\infty} \rtimes_\rho \mathbb{N}$  \cite[section 2]{c}. We have the following fundamental result due to Cuntz \cite{c}.

\begin{theorem}
	For  $n \geq 2$ the Cuntz algebra $\mathcal O_n$ is simple and purely infinite.
\end{theorem}
\subsection{Cuntz-Krieger Algebras}
A more general constructions along the same lines as that  of the Cuntz algebras was suggested by Cuntz and Krieger in \cite{ck}.

\begin{definition}
	Let $A$ be a $n \times n$ matrix with entries  either $0$ or $1$,  and no zero  columns or  rows. Let $\{s_i\}_{i = 1}^n$ be a set of nonzero partial isometries, with mutually orthogonal range projections, acting on a Hilbert space $\mathscr{H}$. The  Cuntz-Krieger algebra $\mathcal O_A$ is then defined as the universal C$^*$-algebra with generators $\{s_i\}_{i = 1}^n$ and additional relations $s_i^*s_i = \sum_{j = 1}^n a_{i,j} s_js_j^*$, for $i,j=1,\cdots, n$.
\end{definition}

In the particular case where $A$ is a matrix with $1$ in all entries, the Cuntz-Krieger relations give back the Cuntz relation $1 = s_i^*s_i = \sum_{j = 1}^n  s_js_j^*$, and we get the Cuntz algebra $\mathcal O_n$. Note that, in this case, the fact that the generators are not just partial isometries, but are isometries, is not changing anything, as by universality, $\mathcal O_A$ only depends on the choice of matrix $A$, and not on the choice of generating partial isometries.

Now let $W$ be as in the case Cuntz algebra, and  $M_A := \{\mu  \in W: s_\mu \neq 0  \}$, and let $\Sigma$ be the set of indices $i \in \{1,\dots,n\}$ for which there are distinct elements $\mu = (\mu_1,\dots,\mu_r)$ and $\nu = (\nu_1,\dots,\nu_s)$ in $M_A$ with $\mu_1 = \nu_1 = \mu_r = \mu_s = i$, while $\mu_k,\nu_l \neq i$, for $1 < k < r$ and $1 < l < s$.
We say that an $n \times n$ matrix $A$  satisfies condition (I) if for all $i \in \{1,\dots, n\}$, there exists an $r-$tuple $(i_1,\dots, i_n) \in M_A$ such that $i_1 = i$ and $i_r \in \Sigma$.
Cuntz and Krieger observed  that for distinct matrices $A$ and $B$ satisfying condition (I), the corresponding Cuntz-Krieger algebras $\mathcal O_A$ and $\mathcal O_B$ are non isomorphic \cite{ck}. It was already observed by Cuntz that the Cuntz algebras $\mathcal O_n$ are non isomorphic for distinct values of $n$, and non of these ``uniqueness'' results are trivial. 

In order to state the simplicity theorem for Cuntz-Krieger algebras, let us recall that a $n \times n$ matrix $A$ is called {\it irreducible} if for all $1 \leq i,j  \leq n$, there is $m \in \mathbb{N}$ such that $(A^m)_{i,j} > 0$, that is the $(i,j)$-th entry of a power of $A$ is non zero \cite{ck}.
 
\begin{theorem}
	Let $A$ be a  $n \times n$  matrix with entries $0$ or $1$, satisfying condition (I), then the Cuntz-Krieger $C^*$-algebra $\mathcal O_A$ is simple if and only if $A$ is irreducible.
\end{theorem}

\subsection{Graph Algebras}

A directed graph $E = (E^0,E^1)$ consists of a vertex set $E^0$, an edge set $E^1$, and range and source maps $r,s : E^1 \rightarrow E^0$. One can define the graph C$^*$-algebra $C^*(E)$ of  $E$, again as a universal C$^*$-algebra. A graph is called row-finite if every vertex receives only finitely many edges, in other words, for each $v \in E^0$, $|r^{-1}(v)| < \infty$. A vertex which receives no edge is called a source. A path is a (finite or infinite) sequence of edges $e_1,e_2,\dots$ with $s(e_{i+1}) = r(e_i)$, for $i = 1,\dots, n-1$. For a path $\mu = e_1,\dots ,e_n$, set   $s(\mu)= s(e_1)$ and $r(\mu) = r(e_n)$. By a cycle we mean a path $\mu$ with $s(\mu) = r(\mu)$.

\begin{definition}
Let $E = (E^0,E^1)$ be a row-finite directed graph and $\mathscr{H}$ be a Hilbert space. A Cuntz-Krieger family $\{S,P\}$ consists of a set of partial isometries $S = \{ s_e \}_{e \in E^1}$ with mutually orthogonal range projections and a set of mutually orthogonal projections $P = \{ p_v \}_{v \in E^0}$, satisfying  relations, $p_v = \sum_{e \in r^{-1}(v)} s_e s_e^*,$ and  $s_e^* s_e = p_{s(e)}$, provided that  $v$ is not source. The 
graph C$^*$-algebra $C^*(E)$ is the corresponding universal $C^*$-algebra with respect to a Cuntz-Krieger family $\{S,P\}$.  
\end{definition}

For every row-finite graph $E$, there is a universal C$^*$-algebra C$^*(E)$ generated by a given Cuntz-Krieger family $\{S,P\}$, in the sense that for any C$^*$-algebra $\mathcal B$ containing a Cuntz-Krieger family $\{ S',T' \}$, there exist a $*$-homomorphism $\pi : C^*(E) \rightarrow \mathcal B$, mapping $S$ and $P$ to $S'$ and $T'$, respectively  \cite[Proposition 1.21]{r}. 

As for the ``uniqueness'', 
	let $E$ be a row-finite directed graph such that every cycle has an entry point, and $\{ S,P\}$ be a Cuntz-Krieger family in a C$^*$-algebra $\mathcal B$ with all projections in $P$ non-zero, then there is an isomorphism between $C^*(E)$ and the C$^*$-algebra $C^*(S,P)$ generated by $S\cup P$ in $\mathcal B$.

To see how this class of algebras extends the class of Cuntz-Kerieger algebras, take a natural number $n \geq 2$, and a matrix $A$ with  entries $0$ and $1$, and define $E_A^0 = \{ 1,2,\dots , n \}$ and $E_A^1 = \{ \overline{ij}:\ a_{i,j} = 1  \}$ with $s(\overline{ij}) = j$ and $r(\overline{ij}) = i$. Let $S = \{ s_i \}_{i = 1}^n$ be a set partial isometries satisfying $s_i^*s_i = \sum_{j = 1}^n a_{i,j} s_js_j^*$ and let $P = \{ s_i s_i^* \}_{i = 1}^n$, which is a set of orthogonal projections, then $\{ S,P \}$ is a Cuntz-Krieger family of  $E_A$ and $\mathcal O_A$ is nothing but $C^*(E_A)$. In particular, one could realize a Cuntz algebra as a graph $C^*$-algebra with the above recipe, but let us give an alternative (but equivalent) way here: for a natural number $n \geq 2$, let  $F = (F^0,F^1)$ be the graph with only one vertex $F^0 = \{ v \}$ and $n$ edges $F^1 = \{ e_1, e_2, \dots , e_n\}$, having trivial range and source maps. Then $\mathcal O_n$ is isomorphic to $C^*(F)$. To illustrate these constructions, let us look at the particular case of $O_2$: the left graph below is the graph just described, whereas the right graph is the graph constructed by the Cuntz-Kriger matrix recipe.
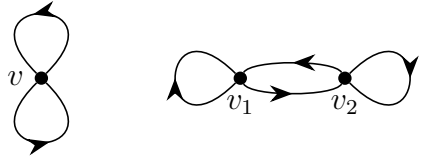
\begin{figure}[h]
	\centering
	\begin{tabular}[r]{@{}cc@{}}
		\begin{tikzpicture}[
		node distance = 12mm and 12mm,
		decoration = {markings, mark=at position 0.55 with {%
				\arrow{Stealth[length=3mm]}},
		},
		dot/.style = {circle, fill, minimum size=5pt,
			inner sep=0pt, outer sep=0pt,
			node contents={}},
		every edge/.style = {draw, semithick,postaction={decorate}},
		]
		
		\node (v) [dot,label= left:$v$, below];
		\path   (v) edge [out=45, in=135, distance=15mm]  (v);
		\path   (v) edge [out=225, in=-45, distance=15mm] (v);
		\end{tikzpicture}
		&
		\begin{tikzpicture}[
		node distance = 12mm and 12mm,
		decoration = {markings, mark=at position 0.5 with {%
				\arrow{Stealth[length=3mm]}},
		},
		dot/.style = {circle, fill, minimum size=5pt,
			inner sep=0pt, outer sep=0pt,
			node contents={}},
		every edge/.style = {draw, semithick,postaction={decorate}},
		]
		
		\node (v) [dot,label= below:$v_1$];
		\node (u) [dot,label= below:$v_2$, right=of v];
		\path   (v) edge [out=225, in=135, distance=15mm]  (v);
		\path   (u) edge [out=45, in=-45, distance=15mm]  (u);
		\path   (u) edge [out=120, in=60, distance=2mm]  (v);
		\path   (v) edge [out=-60, in=240, distance=2mm]  (u);
		\end{tikzpicture}
	\end{tabular}
\end{figure}

Note that the right graph is nothing but the dual graph of the left one, and indeed the C$^*$-algebra of the dual graph is always isomorphic to that of the original graph. More precisely, 
	let $E$ be a row-finite graph with no sources and let $\hat{E}$ be the dual graph with $\hat{E}^0 := E^1$, and $\hat{E}^1$:= the paths of length 2 (with range $r(e_1e_2) = e_1$ and $s(e_1e_2) = e_2$). Then the dual graph is row-finite and gives an isomorphic graph $C^*$-algebra. 
	
Before stating  the simplicity result for graph algebras  \cite[Theorem 4.14]{r}, we need to recall the concept of cofinality. Let $u$ and $v$ be two vertices of a directed graph. By writing $u \leq v$ we mean that there is a path $\mu$ from $v$ to $u$, that is, $s(\mu) = v$ and $r(\mu) = u$. We say that a graph $E$ is {\it cofinal} if for every path $\mu$ and vertex $u \in E^0$, there exists a vertex $v$ on $\mu$ with $u \leq v$. 

\begin{theorem}
	Let $E$ be a row-finite directed graph, then $C^*(E)$ is simple if and only if $E$ is cofinal and every cycle has an entry.
\end{theorem}

In the context of simple graph C$^*$-algebras, there is an important dichotomy: let $E$ be a row-finite graph and $C^*(E)$ is simple, then either $C^*(E)$ is an AF-algebra (if $E$ has no cycle) or it is purely infinite (if $E$ has cycle). Here, being an AF-algebra means that $C^*(E)$ is a direct limit of finite dimensional $C^*$-algebras. 

\subsection{Cuntz-Pimsner algebras and Katsura construction} \label{subsection-katsura}
Let $\mathcal A$ be a C$^*$-algebra and $X$ a right Hilbert $C^*$-module over $\mathcal A$ with right inner product $\langle\cdot,\cdot\rangle_X$. Let $\LL(X)$ be the set of adjointable operators on $X$, then $X$ together with a left action given by a $*$-homomorphism $\phi_X : \mathcal A \rightarrow \LL(X)$ is called a  $C^*$-correspondence over $\mathcal A$. Every $C^*$-algebra is a $C^*$-correspondence over itself by trivial product and inner product $\langle a,b\rangle_\mathcal A := a^*b$. We say that $X$ is  a Hilbert $A$-bimodule if there is a left inner product $_X\la \cdot , \cdot \ra : X \times X \rightarrow \mathcal A$, satisfying $_X\la x , x \ra \geq 0$, and $$_X\la a\cdot x , y \ra = a _X\la x , y \ra,\ _X\la x , y \ra^* =\la y , x\ra_X, \ _X\la x, y \ra z = x\la y,z \ra_X,$$
for each $x,y,z \in X$, $a \in\mathcal  A$. A morphism of $\mathcal A$-correspondences is a bimodule map $T: X_1 \rightarrow X_2$ satisfying $_{X_2}\la T(x), T(y) \ra  =$$_{X_1}$$\la x , y \ra$, for  $a \in\mathcal  A$ and $x,y \in X_1 $, which is the isomorphism of the category of correspondences when $T$ is also surjective, writing $X_1\simeq X_2$. If $T$ is adjointable, and $T^*T = id_{X_1}$, $TT^* = id_{X_2}$, then we have a unitary equivalence, writing $X_1 \approx X_2$.

 We denote by $\K(X)$ the subalgebra of  compact operators on $X$. This is the closed subalgebra generated by ``finite rank'' operators (with the warning that their ranges are merely finitely generated, and not necessarily finite dimensional; in particular, operators in $\K(X)$ are not necessarily compact as operators on the Banach space $X$).  A C$^*$-correspondence $X$ over $\mathcal A$ is called essential  if $\{\phi_X(a)x:\ a \in A, x \in X \}$ is a total set in $X$.

	Let $X$ be a $C^*$-correspondence over $\mathcal A$. By a representation $(\pi,t)$ of $X$ in a $C^*$-algebra $\mathcal B$ we mean a $*$-homomorphism  $\pi :\mathcal  A \rightarrow \mathcal B$ and a linear map $t : X \rightarrow \mathcal  B$ satisfying,
$t(\xi)^*t(\eta) = \pi(\la\xi,\eta\ra)$ and $\pi(a)^*t(\xi) = t(a\cdot\xi)$, for $a\in \mathcal A$ and $\xi,\eta\in X$.  
	The $^*$-algebra $C^*(\pi,t)$ is then generated by the ranges of the maps $\pi$ and $t$.
For a representation $(\pi,t)$ of $X$ in $\mathcal B$, one can define $\psi_t : \mathscr K (X) \rightarrow \mathcal B$, mapping each rank one operator $\theta_{\xi,\eta}$ to $t(\xi)t(\eta)^*$.
The linear maps $\psi^{(\pi,t)}_n : X^{\otimes n} \rightarrow \mathcal B$  are then mapping $x_1 \otimes x_2 \otimes \dots \otimes x_n$ to $t(x_1)t(x_2) \dots t(x_n)$ \cite[section 4]{k}. 
 Now for the Katsura ideal,
\[
	J_X = \phi_X^{-1}(\mathscr K(X)) \cap (ker \phi_X)^\perp,
\]
a Cuntz-Pimsner covariant representation is one for which $\pi(a) = \psi_t(\phi_X(a))$, for all $a \in J_X$.
\begin{definition}
	Let $X$ be a $C^*$-correspondence over $C^*$-algebra $\mathcal A$. The Cuntz-Pimsner algebra $O_X$ is defined as the universal $^C*$-algebra given by Cuntz-Pimsner covariant representations. 
\end{definition}
Katsura in \cite{k} proved a gauge invariance uniqueness theorem and built universal covariant representation $(\pi_X,t_X)$, with $\mathcal O_X \cong  C^*(\pi_X,t_X)$, in the sense that, for all covariant representations  $(\pi,t)$ of $X$, there exist a homomorphism $\rho :\mathcal  O_X \rightarrow  C^*(\pi,t)$ with  $\pi = \rho \circ \pi_X$ and $t = \rho \circ  t_X$.

The graph algebras (and so the Cuntz and Cuntz-Krieger algebras) are  special cases of this construction. To see this, let $E = (E^0,E^1)$ be a row-finite directed graph. Let $a \in \mathcal A:= c_0(E^0)$, $f, g \in c_c(E^1)$, $v \in E^0$ and $e \in E^1$ be given, and define,
\begin{align*}
	(f\cdot a)(e) := f(e)a(s(e)) && \la f,g \ra_\mathcal A (v) = \sum_{e \in s^{-1}(v)} \overline{f(e)}g(e)
\end{align*}
The completion $X_E$ of $c_c(E^1)$ is then a Hilbert $\mathcal A$-module. Defining homomorphism $\phi :\mathcal  A \rightarrow \mathscr L (X_E)$ by $\phi(a)(f)(e) = (a\cdot f)(e) := a(r(e))f(e)$, for $a \in \mathcal A$, $f \in c_c(E^1)$ and $e \in E^1$, we may equip  $X_E$ with a $C^*$-correspondence structure,   referred to as the graph correspondence of $E$. Then  $\phi^{-1}(\K(X_E))$ is the span of point mass functions $\delta_v$ with $| r^{-1}(v)| <\infty$, and $\phi(\delta_v) = 0$, for source vertex $v$ \cite[Proposition 4.4]{fr}. 

Indeed, for a vertex $v$ with $|r^{-1}(v)| < \infty$, 
$
	\phi(\delta_v) = \sum_{e \in r^{-1}(v)} \theta_{\delta_e,\delta_e}.
$
 In this case, the Katsura ideal is nothing but the closed linear span of the set $\{ \delta_v:\ 0 < |r^{-1}(v)| < \infty \}.$
In particular,  a representation $(\pi , t)$ is covariant if and only if for every vertex $v$ with $|r^{-1}(v)| < \infty$, we have,
\[
	\pi(\delta_v) = \psi_t(\phi(\delta_v)) = \sum_{e \in r^{-1}(v)} t(\delta_e) t(\delta_e)^*.
\]
Now setting $S = \{t(\delta_e)\}_{e \in E^1}$ and $P = \{ \pi(\delta_v) \}_{v \in E^0}$, the above computation shows that every covariant representation of $X_E$ induces an Cuntz-Krieger family, and by universality, $C^*(E)$ is nothing but $\mathcal  O_{X_E}$.

For stating simplicity result for Katsura algebras, we need some preparation. Let $X$ be a $C^*$-correspondence over a $C^*$-algebra $\mathcal A$. Then, $X$ is called {\it minimal} if $\la X , J X \ra \subseteq J$, for a closed ideal $J$ of $A$, implies $J = 0$ or $A$.
		It  is called full if  $\{ \la x,y \ra: x,y \in X \}$ is a total set in $\mathcal A$. When $\mathcal A$ is unital, $X$ is called {\it  aperiodic} if for every non zero integer $n$, $X^{\otimes n} \approx I_\mathcal A$ implies $n = 0$, where $X^{\otimes n}$ is the $n$-fold tensor product, and $I_\mathcal A$ is the identity $C^*$-correspondence, that is $\mathcal A$ itself with its canonical inner product.

The first simplicity result  \cite[Theorem 3.9]{s}, is for Cuntz-Pimsner algebras, i.e., Katsura algebras of a full $C^*$-correspondence with injective left action.
 
\begin{theorem}
	Let $X$ be a full $C^*$-correspondence $X$ over a unital $C^*$-algebra $\mathcal A$ with an injective left $\mathcal A$-action. Then $\mathcal O_X$ is simple if and only if $X$ is aperiodic and minimal.
\end{theorem}

 Ery\"{u}zl\"{u} and Tomforde in \cite{et} defined the condition $(S)$ (where $S$ stands for simplicity) and obtained a  simplicity result for Katsura  algebras  \cite[Theorem 4.3]{et}.
	Let $X$ be a $C^*$-correspondence over a $C^*$-algebra $\mathcal A$. Let $(\pi_X,t_X)$ be a universal covariant representation of $X$ into $\mathcal O_X$. A vector $\xi \in x^{\otimes n}$ with $n \in \mathbb{N}$ is called non-returning if for $m < n$ and $\eta \in X^{\otimes m}$, $t^{\otimes n}_X(\xi) t^{\otimes m}_X(\eta) t_X^{\otimes n}(\xi) = 0$. We say $X$ satisfies condition $(S)$ if for each $a \in \mathcal A$ with $a \geq 0$, and each $n \in \mathbb{N}$ and $\epsilon > 0$, there exist $m > n$ and a non-returning unit vector $\xi \in X^{\otimes m}$ with $\| \la \xi , a\xi \ra \| > \| a \| - \epsilon$. An ideal $I$ of $\mathcal A$ is called invariant  if $IX \subseteq XI$, and for all $a \in J_X$, $aX \subseteq XI$ implies $a \in I$.

\begin{theorem}
	Let $X$ be a $C^*$-correspondence over a C$^*$-algebra $\mathcal A$. If $X$ satisfies condition $(S)$ and $\mathcal  A$ has no non trivial invariant ideal, then $\mathcal O_X$ is simple.
\end{theorem}

\subsection{Higher rank graph algebras}
As a generalization of the graph algebras, Kumjian and Pask introduced the notion of Higher rank graphs in \cite{kp}.
	For a natural number $k$, a $k$-graph (or a graph of rank $k$) $(\Lambda,d)$ consist of a countable small category $\Lambda$, with source and range maps $s,r$, and a functor $d : \Lambda \rightarrow \mathbb{N}^k$ satisfying the factorization property, that is, for all $\lambda \in \Lambda$, $m,n \in \mathbb{N}^k$ with $d(\lambda) = m+n$, there exist unique elements $\mu,\nu \in \Lambda$ such that $\lambda = \mu \nu$ with $d(\mu) = m$ and $d(\nu) = n$.

Often for $n \in \mathbb{N}^k$, the notation $\Lambda^n$ is used for $d^{-1}(n)$. Similar to the directed graph case, a $k$-graph is called row-finite if $r^{-1}(v) \cap \Lambda^n$ is finite for every $v \in \Lambda^0$ and $n \in \mathbb{N}^k$. 

As a concrete example, for a natural number $k$, define $\Omega_k$ to be the small category with object set $\mathbb{N}^k$ and morphism set consisting of $(m,n) \in \mathbb{N}^k \times \mathbb{N}^k$ with $m \leq n$, with the range and source maps  defined by $r(m,n) = m$ and $s(m,n) = n$. For $d : \Omega_k \rightarrow \mathbb{N}^k;\ d(m,n) = n-m$, $(\Omega_k, d)$ is a $k$-graph.

In order to build a universal $C^*$-algebra like the previous cases,  we need to define proper relations.
	For a row-finite $k$-graph $(\Lambda,d)$, a Cuntz-Krieger-$\Lambda$ family is a family $S = \{s_\lambda: \ \lambda \in \Lambda\}$ of partial isometries, with $\{ s_v:\ v \in \Lambda^{0} \}$ being  mutually orthogonal projections, satisfying, 
		 $s_\lambda s_\mu = s_{ \lambda \mu}$ for  $s(\lambda) = r(\mu)$, 
		 $s_\lambda^* s_\lambda = s_{s(\lambda)}$, and 
		 $s_v = \sum_{\lambda \in r^{-1}(v) \cap \Lambda^n} s_\lambda s_\lambda^*$, for  $v \in \Lambda^0$ and $n \in \mathbb{N}^k$.

The higher rank graph algebra $C^*(\Lambda)$ is then the universal $C^*$-algebra generated by Cuntz-Krieger-$\Lambda$ families. Every row-finite $k$-graph $\Lambda$ with no source, is generated by a universal Cuntz-Krieger-$\Lambda$ family $\{s_\lambda\}_{\lambda \in \Lambda}$, in the sense that, for any Cuntz-Krieger-$\Lambda$ family $\{t_\lambda\}_{\lambda \in \Lambda}$ in a $C^*$-algebra $\mathcal B$, there exist a homomorphism $\pi : C^*(\Lambda) \rightarrow \mathcal B$ with $\pi(s_\lambda) = t_\lambda$.

To see that a row-finite directed graph is a special case, we need to build a small category $\mathscr{P}_E$ (path category) out of a directed graph $E = (E^0,E^1)$ with source and range maps $s_E$ and $r_E$. Let the object set $\mathscr{P}_E^0$ of $\mathscr{P}_E$ consist of  vertices in $E^0$ and the morphisms  of $\mathscr{P}_E$ simply be the finite paths $\mu$ with $s(\mu) = s_E(\mu)$ and $r(\mu) = r_E(\mu)$. Consider the product  $(\mu,\nu) \rightarrow \mu\nu$, defined when $s(\mu) = r(\nu)$, by pasting the paths. Define the functor of $d : \mathscr{P}_E \rightarrow \mathbb{N}$ to be the path length function. The directed graphs are then just $1$-graphs in this context.

The simplicity of higher rank graph C$^*$-algebras was characterized  by Kumjian and Pask \cite[Proposition 4.8]{kp}. First, we need to recall some definitions.
	Let $\Lambda$ be a $k$-graph, and define the infinite path space of $\Lambda$ as the set $\Lambda^\infty$ of all $ k$-graph morphisms $x: \Omega_k \rightarrow \Lambda$. 
		For  $v \in \Lambda^0$ put $\Lambda^\infty(v) = \{ x \in \Lambda^\infty: x(0) = v \}$, and for $p \in \mathbb{N}^k$, define $\sigma^p : \Lambda^\infty \rightarrow \lambda^\infty$ by $\sigma^p(x)(m,n) = x(m+p, n+p)$, for $x \in \Lambda ^\infty$ and $(m,n) \in \Omega_k$.
		 We say that $\Lambda$ is {\it cofinal} if for every $x \in \Lambda^\infty$ and vertex $v \in \Lambda^0$, there is $\lambda \in \Lambda$ and $n \in \mathbb{N}^k$ with $s(\lambda) = x(n)$ and $r(\lambda) = v$.
		 A path $x \in \Lambda^\infty$ is periodic if there exist $p \in \mathbb{Z}^k$ such that for  $(m,n) \in \mathbb{N}^k$ with $m+p \geq 0$, $x(m+p, n+p) = x(m,n)$; and aperiodic if $\sigma^nx$ is not periodic for all $n \in \mathbb{N}^k$. We say that $\Lambda$ satisfies the aperiodicity condition, if for every vertex $v \in \Lambda^0$ there is an aperiodic path $x \in \Lambda^\infty(v)$.

\begin{theorem}
	If $\Lambda$ satisfies the aperiodicity condition, then $C^*(\Lambda)$ is simple if and only if $\Lambda$ is cofinal.
\end{theorem}

\section{Product system and Cuntz-Nica-Pimsner algebras}
 
Let us review the notion of product system and construction of the full and reduced Cuntz-Nica-Pimsner algebras $\mathcal{NO}_X$ and $\mathcal{NO}^r_X$. The $C^*$-algebras of higher rank graphs and Cuntz-Pimsner algebras fit in this framework.

 	Let $\mathcal A$ be a $C^*$-algebra, $P$ be a discrete monoid with identity $e$,  and for each $p\in P$, $X_p$ be a $C^*$-correspondence over $\mathcal A$. A product system of $C^*$-correspondences over $P$ is a semigroup (under tensor product) $X = \{ X_p \}_{p \in P}$ consisting of $C^*$-correspondences, with $X_e\simeq A$, as $C^*$-correspondences,  such that the multiplication in $X$ by elements of $X_e$ implements the left and right actions of $\mathcal A$ on each $X_p$, and for $p,q \in P \backslash \{e\}$,
 	there is an isomorphism of $\mathcal A$-correspondences $M_{p,q} : X_p \otimes_\mathcal A X_q \rightarrow X_{pq}$,  satisfying $M_{p,q}(x \otimes_A y) = xy$, for  $x \in X_p$ and $y \in X_q$. 
 When $ p \neq e$, there is a morphism $\iota_p^{pq} : \LL(X_p) \rightarrow \LL(X_{pq})$ satisfying $
 	\iota_p^{pq}(S)(xy) = (Sx)y $, for  $x \in X_p, y \in X_q$, and $S \in \LL(X_p)$.  
 If we identify $\K(X_e)$ with $\mathcal A$, we also get $\iota_e^p : \K(X_e) \rightarrow \LL(X_q)$, given by the left action $\phi_p$ of $X_p$.
 
In this review, we  restrict ourselves to monoids arising from quasi-lattice ordered groups. Let $G$ be a discrete group with identity $e$, $P \leq G$ be a sub-semigroup of $G$ with $P \cap P^{-1} = \{e\}$. The pair $(G,P)$ is called a quasi-lattice ordered group if for the partial order on $G$ defined for  $p,q \in G$ by $p\leq q$ if $p^{-1}q \in P$, any two elements $p,q\in P$ with a common upper bound have a least common upper bound $p \lor q \in P$. We write $p\lor q = \infty$ when $p$ and $q$ fail to have a common upper bound in $P$ and $p\lor q < \infty$, otherwise.  We denote  the $A$-valued inner product on $X_p$ by $\la \cdot,\cdot \ra_p$. 

A representation $\psi$ of $X$ in a $C^*$-algebra $ \mathcal B$ is a map $\psi: X \rightarrow \mathcal B$ such that $\psi_e$ is a $*$-homomorphism, and for $p \in P$, $\psi_p := \psi|_{X_p}$ is linear, satisfying $\psi_p(x)\psi_q(z) = \psi_{pq}(xz)$ and $\psi_e(\la x, y \ra_p) = \psi_p(x)^*\psi_p(y)$, for  $p,q \in P$, $x,y \in X_p$ and $z \in X_q$. A representation $\psi$ is called injective if  $\psi_e$ is injective. In this case, it follows that each $\psi_p$ is an isometry.
For each $p\in P$, there is a $*$-homomorphism $\psi^{(p)}: \K(X_p) \rightarrow \mathcal B$ satisfying, $	\psi^{(p)}(\theta_{x,y}) = \psi_p(x)\psi_p(y)^*$, for $x,y \in X_p$ \cite{pim}. When the product system is a Hilbert C$^*$-bimodule, a representation of $X$ is said to be a Hilbert $C^*$-bimodule representation of $X$ if moreover, $
	\psi_e(_{p}\la x, y \ra) = \psi_p(x)\psi_p(y)^*$,  for $p \in P$ and $x,y \in X_p$. 
The system $X$ 
is called compactly aligned if $\iota_p^{p\lor q}(S)\iota_q^{p\lor q}(T) \in \K(X_{p\lor q})$, for each  $p \lor q < \infty$, $S \in \K(X_p)$,  and $T \in \K(X_q)$.
For the rest of the section, $(G,P)$ is a quasi-lattice ordered group and $X$ is a product system over $P$.
 
We denote  the $C^*$-algebra generated by the image of a representation $\psi$ in $\mathcal B$ by $C^*(\psi)$. If $\iota$ is the universal representation of product system $X$ over $P$ and  $\mathcal{T}_X := C^*(\iota)$, then for each representation $\psi$ of $X$, there is a $*$-homomorphism $\rho: \mathcal{T}_X \rightarrow C^*(\psi)$ with $\rho(i_p(x)) = \psi_p(x)$,   $p \in P$ and $x \in X_p$.

	Let $X$ be a compactly aligned product system over $P$. A representation $\psi$ of $X$ is Nica coveariant if for each $p,q\in P$ and $S \in \K(X_p)$ and $T \in \K(X_q)$,
$$		\psi^{(p)}(S)\psi^{(q)}(T) = 
		\begin{cases}
		\psi^{(p \lor q)}(\iota_p^{(p\lor q)}(S)\iota_q^{(p\lor q)}(T)) & \text{ if } p \lor q < \infty
		\\
		0 &\  \text{if } p\lor q=\infty
		\end{cases}$$
 The $C^*$-algebra generated by the universal Nica covariant representation $\iota=\iota_X$ is denoted here by $\mathcal T^{\rm cov}_X$ (or by $\mathcal{T}_{\rm cov}(X)$ in some texts). It is known that the set $\{\iota(x)\iota(y)^*:  x,y \in X\}$ is total in 
$	\mathcal{T}^{\rm cov}_X$ \cite[2.2]{clsv}.

\subsection{Full and reduced Cuntz-Nica-Pimsner algebras}
Let $(G,P)$ be a quasi lattice ordered group and $X$ is a compactly aligned product system over $P$. For the closed ideal ideals $I_e = \mathcal A$ and $I_q = \bigcap_{e < p \leq q} ker(\phi_p)$, for $q \in P \backslash \{e\}$,  consider the $\mathcal A$-correspondence  $
\tilde{X}_q = \bigoplus_{p\leq q} X_p . I_{p^{-1}q}$, 
with implementing left action $\tilde{\phi}_q$. We say that $X$ is $\tilde{\phi}$-injective if all the left actions $\tilde{\phi}_q$ are injective. For  $p,q \in P$, let  $\tilde{\iota}_p^q: \LL(X_p) \rightarrow \LL(\tilde{X}_q)$ be defined by
$
	\tilde{\iota}_p^q(T) = \bigoplus_{r \leq q} \iota_p^r(T).
$

A representation $\psi$ of $X$ is called Cuntz-Nica-Pimsner covariant (CNP-covariant) if for each finite subset $F \subseteq P$, each $p \in F$ and $T_p \in \K(X_p)$, if $\sum_{p \in F} \tilde{\iota}_p^q(T_p) = 0$, then $
		\sum_{p \in F} \psi^{(p)}(T_p) = 0, $  for large $q$, that is, for each $s \in P$ there exists $r \geq s$ such that the latter sum is zero for $r \leq q$.

The $C^*$ algebra $\mathcal{NO}_X$ generated by the range of universal CNP-covariant representation $j_X$ is called the (full) Cuntz-Nica-Pimsner algebra of $X$. Since $j_X$ is an injective representation, by universality of $\mathcal{T}^{\rm cov}_X$, there exists a canonical surjective $*$-homomorphism $\rho_{CNP}: \mathcal{T}^{\rm cov}_X \rightarrow \mathcal{NO}_X$. 
Let $\iota=\iota_X$ be the universal Nica covarinat representation generating $\mathcal{T}^{\rm cov}_X$, then the core $\mathcal F_X$ of $\mathcal{T}^{\rm cov}_X$ is the span closure $\{ i_p(x)i_p(y)^* | x,y \in X_p \}$, which also coincides with the span closure of ranges of $i_X^{(p)}$. 

A criterion for injectivity of induced representations is available as follows \cite[Theorem 3.8]{clsv}: 
	assume that either  the left actions on all fibers of $X$ are injective or that $P$ is directed and $X$ is $\tilde{\phi}$-injective. For a CNP-covariant representation $\psi$ of $X$ on a $C^*$-algebra $B$ where $\rho: \mathcal{NO}_X \rightarrow B$ is the induced homomorphism, $\rho$ is in injective on $\rho_{CNP}(\mathcal F_X)$ if and only if $\psi$ is injective as a representation.

Let $\delta_G: C^*(G) \rightarrow C^*(G) \otimes C^*(G)$ be the canonical homomorphism  mapping  $g \in G$ to $i_G(g) \otimes i_G(g)$, where $i_G$ is the canonical inclusion of $G$ in the full group $C^*$-algebra $C^*(G)$.	A (full) coaction of $G$ on a $C^*$-algebra $A$ is a nondegenerate injective $*$-homomorphism $\delta: A \rightarrow A \otimes C^*(G)$ satisfying, $
		(\delta \otimes id_{C^*(G)}) \circ \delta = (id_A \otimes \delta_G) \circ \delta.$ The generalized fixed-point algebra of $\delta$ is $A_e^\delta$, where  $A_g^\delta = \{a \in A | \delta(a) = a \otimes i_G(g)  \}$, for $g \in G$. The generalized fixed-point algebra associated with $\mathcal{T}^{\rm cov}_X$ is nothing but the core of $\mathcal{T}^{\rm cov}_X$ \cite{clsv}.
A representation $\psi$ of $X$ in $A$ is said to be  gauge compatible if there is a coaction $\delta$ of $G$ on $A$ with $
\delta(\psi_p(x_p)) = \psi_p(x_p) \otimes i_G(p)$, for  $p \in P$ and $x \in X_p$. The $C^*$-algebras $\mathcal{T}^{\rm cov}_X$ and $\mathcal{NO}_X$ associated to a compactly aligned product system $X$ over $P$ have gauge compatible coactions with respect to their universal representation \cite{clsv}. 

Again, let us suppose that either the left action on each fiber is injective or  $P$ is directed. Also assume that $X$ is $\tilde{\phi}$-injective. Then the $C^*$-algebra $\mathcal{NO}_X^r$ generated by the couniversal, gauge compatible injective Nica covariant representation of $X$  is called the reduced Cuntz-Nica-Pimsner $C^*$-algebra of $X$, whose canonical coaction is denoted by $\nu^r$. The $C^*$-algebras $\mathcal{NO}_X$ and $\mathcal{NO}^r_X$ are known to be the same as the full and reduced cross sectional algebras  of the Fell bundle $\{ (\mathcal{NO}_X)_g^\nu \times \{g\} \}_{g\in G}$ \cite{clsv}. When $X$ is  $\tilde{\phi}$-injective,  we say that $\mathcal{NO}_X$ has the gauge invariant uniqueness property if for any $C^*$-algebra $B$, the injectivity of any surjective $*$-homomorphism $\phi: \mathcal{NO}_X \rightarrow B$ is equivalent to  the injectivity of $\phi|_{\iota_X(A)}$ plus the existence of  a coaction $\beta$ of $G$ on $B$ with $\beta \circ \phi = (\phi \otimes \text{id}_{C^*(G)}) \circ \nu$. This is known to be automatic whenever $G$ is amenable \cite[Corollary 4.12]{clsv}. In general, the following guage invariance uniqueness property is verified in \cite{clsv}: suppose  that either the left action on each fiber is injective or $P$ is directed, and $X$ is $\tilde{\phi}$-injective, then the following conditions are equivalent:

$(i)$ the canonical surjection of $\mathcal{NO}_X$ onto $\mathcal{NO}^r_X$ is an isomorphism,
		
$(ii)$ $\mathcal{NO}_X$ has the gauge invariant uniqueness property,

$(iii)$ given two injective gauge compatible CNP-covariant representations  $\psi_i : X \rightarrow B_i$ of $X$ whose image generates $B_i$, $i=1,2$, there is a $*$-isomorphism $\phi : B_1 \rightarrow B_2$ with $\phi \circ \psi_1 = \psi_2$.

As promised, we now explain how Cuntz-Pimsner algebras and higher rank graph algebras could be considered as Cuntz-Nica-Pimsner algebras. Let $X$ be a $C^*$-correspondence over a $C^*$-algebra $\mathcal A$, then the product system $X^\otimes$ over $\mathbb{N}$ could be defined as the monoid of $n$-fold tensor products $X^{\otimes n}$ with $X^{\otimes 0} :=\mathcal  A$, where the isomorphism $M_{m,n} : X^{\otimes m} \otimes X^{\otimes n} \rightarrow X^{\otimes m+n}$ gives the product of the monoid. Since $\mathbb{N}$ is totally ordered, the product system is compactly aligned. Now we could describe how Sims and Yeend proved that the Cuntz-Pimsner algebras are special cases of Cuntz-Nica-Pimsner algebras \cite[Proposition 5.3]{sy}. Let $X$ be a $C^*$-correspondence over a $C^*$-algebra $\mathcal A$, let $(\pi_X,t_X)$ be the universal covatiant representation of $X$ on $\mathcal O_X$ and $j_X$ be the universal representation of $X^{\otimes}$ on $\mathcal{NO}_{X^{\otimes}}$, then there is an isomorphism $\theta : O_X \rightarrow \mathcal{NO}_{X^{\otimes}}$,  mapping $\pi_X(a)$ and $t_X(x)$ to $j_X(a)$ and $j_X(x)$, respectively. For a representation $(\pi,t)$ of $X$ the representation $\psi^{(\pi,t)}$ (defined as in the subsection \ref{subsection-katsura}) is CNP-covariant if and only if $(\pi,t)$ is covariant.

As a concrete example, let $(\Lambda,d)$ be a row-finite $k$-graph and $A = c_0(\Lambda_0)$. For $n \in \mathbb N^k$, let $X_{\Lambda^n}$ be  defined as $X_E$ (see, subsection \ref{subsection-katsura}). For $\lambda \in \Lambda^n$ and $\mu \in \Lambda^m$, define $\delta_\lambda \delta_\mu=\delta_{\lambda \mu}$, if $(\lambda,\mu)$ is multiplicable, and 0 otherwise. Then $P_\Lambda = \{ X_{\Lambda^n}  \}_{n \in \mathbb{N}^k}$ is a product system over $\mathbb{N}^k$ \cite[Section 5]{sy}). For $\lambda , \nu$, consider the collection $MCE(\mu,\nu)$ consisting of those $\lambda \in \Lambda$ satisfying $d(\lambda) = d(\mu) \wedge d(\nu), \lambda = \mu\mu' = \nu\nu'$, for some $\mu'\nu' \in \Lambda \}$. The $k$-graph $\Lambda $ is then called finitely aligned if $MCE(\mu,\nu)$ is finite, for every choice $\mu , \nu \in \Lambda$. Sims and Yeend \cite[Proposition 5.4]{sy} observed that $C^*$-algebras of higher rank graphs could be described as Cuntz-Nica-Pimsner algebras as follows. 
		Let $\Lambda$ be a finitely aligned $k$-graph, and $P_\Lambda$ be the associated product system. Let  $\{s_\lambda\}_{\lambda \in \Lambda}$ be the universal Cuntz-Krieger-$\Lambda$ family in $C^*(\Lambda)$, and let $j_P$ be the universal $CNP$-covariant representation of $X$ in $\mathcal{NO}_P$. Then there is an isomorphism $C^*(\Lambda) \rightarrow \mathcal{NO}_P$,  mapping $s_\lambda$ to $j_P(\delta_\lambda)$, for all $\lambda \in \Lambda$.

\subsection{Kishimoto Condition and Simplicity}

Throughout this subsection, we assume that either the left action on each fiber is injective or  $P$ is directed, and $X$ is $\tilde{\phi}$-injective. This gurantees the existence of the  reduced Cuntz-Nica-Pimsner $C^*$-algebra   $\mathcal{NO}_X^r$ of $X$ generated by the couniversal, gauge compatible injective Nica covariant representation of $X$. 
As a simple observation, note that for a closed ideal $I$ of $\mathcal{NO}_X$,  $I \cap j_e(A) = \{0\}$ implies $I \cap (\mathcal{NO}_X)^{\delta}_e = \{0\}$.

Next let us define the notions of minimality and Kishimoto condition for a quasi-lattice ordered group $(G,P)$ and product system $X$ over $P$. This is related to the notion of aperiodicity for $C^*$-correspondences, introduced by Muhly and Solel \cite[Definition 5.1]{ms}, inspired by the results of Kishimoto \cite[Lemma 1.1]{ki} and Olesen
and Pedersen \cite[Theorems 6.6 and 10.4]{op} (see also, \cite{gs}). In a more general setting, aperiodicity of Fell bundles \cite[Definition 4.1]{ks} and their intersection property have been studied in \cite{ks}. It is known that, aperiodicity implies the intersection property. Also, the connection between simplicity of crossed sectional algebras of a Fell bundle with intersection property is described in  \cite[Corollary 3.13]{ks}. In our context, aperiodicity of Fell bundles associated with a product system and the Kishimoto condition are basically the same (cf.,  \cite[Lemma 4.2]{ks}). There is an alternative, more direct way to obtain a simplicity result for Cuntz-Nica-Pimsner algebras, and the rest of this subsection  is devoted to a rather detailed description of such a direct approach, which is essentially an adaptation of the argument for the case of crossed products (cf., \cite{p}).

Let us call a product system $X$ over $P$ {\it minimal} if  $\la X_p , J X_p \ra \subseteq J$, for a closed ideal $J$ of $\mathcal A$ implies $J = 0$ or $\mathcal A$, and {\it essential} if each $X_p$ is essential, for $p\in P$. The system  $X$ is said to satisfy the  {\it Kishimoto condition}, if for any norm one positive element $x \in (\mathcal{NO}^r_X)^\nu_e$,  $\varepsilon > 0$, and finite subsets $F\subseteq G \backslash \{ e \} $, $S\subseteq \bigcup_{g \in F} (\mathcal{NO}^r_X)^\nu_g$, there exists norm one positive element $c \in (\mathcal{NO}^r_X)^\nu_e$ with $
		\| cxc \| > 1-\varepsilon, \ \| c x_g c \| < \varepsilon,$
for each  $g \in F$ and  $x_g \in S$.

Since $\mathcal{NO}^r_X$ is the reduced cross sectional algebra of the Fell bundle $\{(\mathcal{NO}^r_X)_g^\nu \times \{g\}\}_{g\in G}$, by the standard theory of Fell bundles, there is a faithful conditional expectation $\mathbb E: \mathcal{NO}^r_X \rightarrow (\mathcal{NO}^r_X)^\nu_e$,  vanishing on all fibers except the one at $e$ \cite[Proposition 2.12]{e}. Also, $\mathcal{NO}^r_X$ is a graded $C^*$-algebra and $\bigoplus_{g \in G} (\mathcal{NO}^r_X)_g^\nu$ is dense in $\mathcal{NO}^r_X$ \cite[Proposition 3.2]{e}.

We have the following characterization of the Kishimoto condition in terms of finite subsets of $\mathcal A$ under the above assumptions: assume that there exist unitaries $u_g \in (\mathcal{NO}^r_X)_g^\nu \cap ((\mathcal{NO}^r_X)_e^\nu)'$, for  $g \in G$, such that $u_{g^{-1}} = u_g^*$. When $X$ is a product system consisting of $\mathcal A$-bimodules,  the Kishimoto condition is equivalent to  the condition that for each norm one positive element $a \in\mathcal A$, finite subsets $F \subset G\backslash \{1\}$,  $S \subset\mathcal A$, and $\varepsilon > 0$, there exists  a norm one positive element $c \in \mathcal A$ with 
	$\| cxc \| > 1-\varepsilon, \ \| c s c \| < \varepsilon, $ for all $g \in F$ and $s \in S$. 

Next, let us observe that in the particular case of product systems coming from actions and corresponding twisted crossed products, the above  Kishimoto condition is nothing but the classical Kishimoto condition for actions: let $\alpha$ be an action of $G$ on a unital $C^*$-algebra $\mathcal  A$ with unit 1, $\omega$ be a $\mathbb{T}$-valued cocycle on $G$, and let $\bar\omega$ be the conjugate cocycle. Put $X_p =\mathcal  A$, for each $p\in P$, and consider $A$ as an $\mathcal A$-Hilbert bimodule with actions and $\mathcal A$-valued inner products, 
	\[
	\la x,y \ra := x^*y, \ p\la x,y \ra_p := \alpha_{p^{-1}} (xy^*), \ a\cdot x = \alpha_p(a)x, \ x\cdot a = xa,
	\]
	for $x,y \in X_p=\mathcal A$ and $a \in \mathcal A$. Define isomorphism $X_p \otimes_\mathcal A X_q \rightarrow X_{pq}$; $x\otimes y\mapsto \bar\omega(q,p)\alpha_q(x)y$, then $X = \{X_p \}_{p \in P}$ is a product system consisting of essential $C^*$-bimodules. If $P$ be a directed and $G$ generated by $P$ as a group, then there is a canonical isomorphism between $\mathcal  A \rtimes_{\alpha,\omega}^r G$ and $\mathcal{NO}_X^r$ \cite[Lemma 5.1]{clsv}.
Let us recall the definition of the classical Kishimoto condition for actions. This notion first appeared in \cite{ki}, where Kishimoto showed that if a discrete group $G$ acts on a $C^*$-algebra $\mathcal A$, say by $\alpha$, with no non trivial $G$-invariant closed ideal (i.e., $\alpha$ is minimal), then the reduced crossed product $\mathcal A\rtimes_r G$ is simple whenever the strong Connes spectrum \cite{ki2} of automorphisms $\alpha_g$ is non trivial, for $g\in G\backslash\{e\}$. Kishimoto showed that the above two conditions together are stronger than the following formulation of the Kishimoto condition \cite[Lemma 3.2]{ki}, a name coined by Chris Phillips in \cite{p}, but it is implicit in his proof of simplicity of the reduced crossed product \cite[Theorem 3.1]{ki}, that minimality of $\alpha$ and the Kishimoto condition are good enough to imply simplicity.
	Let $G$ be a discrete group, $\alpha$ be an action of $G$ on a $C^*$-algebra $\mathcal A$. We say that $\alpha$ satisfies the Kishimoto condition if for all norm one positive elements $a \in \mathcal A$, finite subsets $F \subseteq G\backslash \{1\}$,  $S \subseteq \mathcal A$, and  $\varepsilon > 0$, there exists a norm one positive element  $c \in \mathcal A_{+}$ with  $\| cxc \| > 1-\varepsilon, \ \| c b \alpha_g(c) \| < \varepsilon, $ for all $g \in F$ and $b \in S$.

		Let $P$ be directed and generate $G$, $\alpha$ be an action of $G$ on a unital $C^*$-algebra $\mathcal A$ with unit 1, and $\omega$ be a $\mathbb{T}$-valued cocycle on $G$.  Let $X$ be the product system of the twisted crossed product $\mathcal A \rtimes_{\alpha,\omega}^r G$, then $X$ satisfies the Kishimoto condition if and only if $\alpha$ satisfies the Kishimoto condition. 
It is easy to see that the general Kishimoto condition defined in this section coincides with the classical Kishimoto condition for the special case of twisted crossed products. 

Finally we have the following simplicity result for the Cuntz-Nica-Pimsner algebras. 

\begin{theorem}
	Let $(G,P)$ be a quasi-lattice ordered group, $X$ be a product system over $P$ consisting of $C^*$-correspondences over unital $C^*$-algebra $\mathcal  A$. If $X$ is minimal or $\mathcal A$ is simple, $X$ satisfies the Kishimoto condition and $j^r_e(\mathcal  A)$ contains an approximate identity of $\mathcal{NO}^r_X$, then $\mathcal{NO}^r_X$ is simple.
\end{theorem}
	
	This result follows from more general results on simplicity of cross sectional $C^*$-algebras of Fell bundles \cite[Corollary 3.13]{ks}.

 \end{document}